%


\def\un{\cal}

\def\k{\kappa}

\def\cala{{\cal A}}
\def\calb{{\cal B}}

\def\lessdot{\mathrel{\mathord{<}\!\!\raise 0.8 pt\hbox{$\scriptstyle\circ$}}}

%



\baselineskip=0.27in

\def\vd{\Vert\kern -0.2em \hbox{\rm -}}

\centerline{\bf Complete Quotient Boolean Algebras}
\vskip 0.3in
\centerline{Akihiro Kanamori and Saharon Shelah\footnote*{The second author is
 partially supported by the U.S.-Israel Binational Science Foundation. This
 paper is \#390 in his bibliography.}}
\bigskip
For $I$ a proper, countably complete ideal on ${\cal P}(X)$ for some set $X$,
 can the quotient Boolean algebra ${\cal P}(X)/I$ be complete? This question
 was raised by Sikorski [Si] in 1949. By a simple projection argument as for
 measurable cardinals, it can be assumed that $X$ is an uncountable cardinal 
$\kappa$, and that $I$ is a $\kappa$-complete ideal on ${\cal P}(\kappa )$ 
containing all singletons. In this paper we provide consequences from and 
consistency results about completeness. Throughout, $\kappa$ will denote an
 uncountable cardinal, and by an {\it ideal over} $\kappa$ we shall mean a 
proper, $\kappa$-complete ideal on ${\cal P}(\kappa )$ containing all 
singletons.

If $\kappa$ is a measurable cardinal and $I$ a prime ideal over $\kappa$, then
 of course ${\cal P}(\kappa )/I$ is complete, being the two-element Boolean 
algebra. The following theorem shows that completeness
 in itself has strong consistency
 strength:
\bigskip
{\bf Theorem A.} {\it If $\kappa\geq\omega_3$ and  there is an ideal $I$ over 
$\kappa$ such that ${\cal P}(\kappa )/I$ is complete, then there is an inner
 model with a measurable cardinal.}

The restriction $\kappa\geq\omega_3$ is necessary for our proof.

There is a well-known situation in which completeness obtains. An ideal over
 $\kappa$ is $\lambda${\it -saturated iff\/}
 whenever $\{ X_\alpha\ |\ \alpha <\lambda\}
\subseteq {\cal P}(\kappa ) - I$, there are $\alpha <\beta <\lambda$ such 
that $X_\alpha \cap X_\beta\not\in I$. Smith-Tarski [ST] established that if an
 ideal $I$ over $\kappa$ is $\kappa^{+}$-saturated, then ${\cal P}(\kappa )/I$
 is complete. There are several consistency results about the existence of 
such ideals; modifying a sharp one due to the second author, we establish the
 following result separating completeness from saturation. NS denotes
 the ideal 
over $\omega_1$ consisting of the non-stationary subsets; for an ideal $I$ over
 $\kappa$ and $A\subseteq\kappa$, $I|A$
 denotes $\{ X\subseteq\kappa\ |\ X\cap A\in I\}$,
 the restriction of $I$ to $A$; finally, for $S\subseteq\omega_1,\ \tilde S$
 denotes $\omega_1 - S$.
\bigskip
{\bf Theorem B.} {\it Suppose that $\kappa$ is a Woodin cardinal, {\rm CH} 
holds, and
 $S$ is a stationary, co-stationary subset of $\omega_1$. Then there is a 
cardinal-preserving forcing extension with no new reals satisfying: 
 $2^{\omega_1}=\kappa$, and $I={\rm NS}|\tilde S$ is an ideal over 
$\omega_1$ such that ${\cal P}(\omega_1)/I$ is a complete Boolean algebra (and
 $I$ is not $\omega_2$-saturated).} 

Of course, starting with a Woodin cardinal we can insure that CH also holds
 by carrying out a preliminary extension. Also,
 the statement has the same force if $S$ is replaced by $\tilde S$, but the 
 formulation is notationally convenient for the proof. That the $I$ in
 the resulting model cannot be $\omega_2$-saturated follows from the following
 well-known result of Ketonen [Ke]: If CH holds and there is an 
$\omega_2$-saturated ideal over $\omega_1$, then $2^{\omega_1}=\omega_2$.
 Adding a further layer of complexity to the proof of
 Theorem B, we also establish:
\bigskip
{\bf Theorem C.} {\it Assume the hypotheses of Theorem B. Then there is a 
forcing extension with no new reals satisfying: 
$2^{\omega_1}=\omega_3=\kappa$, and $I={\rm NS}|\tilde S$ is an ideal over 
$\omega_1$
 such that ${\cal P}(\omega_1)/I$ is a complete Boolean algebra (and $I$ is
 not $\omega_2$-saturated).}

It follows from 1.1(a) below that $I$ must be $\omega_3$-saturated.
 This then contributes to the theory of saturated ideals by establishing the
 relative consistency of $2^\omega =\omega_1,\ 2^{\omega_1}=\omega_3$, and the
 existence of an $\omega_3$-saturated ideal over $\omega_1$.

In $\S$1 we derive some consequences of completeness and establish Theorem A.
 In $\S$2 we indicate the modifications necessary to a previous proof of the
 second author to establish Theorem B. Finally in $\S$3 we build on $\S$2 to
 establish Theorem C, providing  iteration lemmas for iterated semiproper
 forcing with mixed supports. The main mathematical advances in this paper are
 due to the second author, based on  speculations and prodding by the first.
\vskip 0.3in
\centerline{$\S$1 Consequences of Completeness}
\bigskip
We first review the various concepts involved to affirm some notation: Let $I$
 be an ideal over $\kappa$. Then $I^{+}={\cal P}(\kappa ) - I$, the 
``positive measure'' sets with respect to $I$. For any $A\subseteq\kappa,\ 
[A]_I=\{ B\subseteq\kappa\ |\ A\bigtriangleup B\in I\}$. The Boolean algebra ${\cal P}
(\kappa )/I$ consists of the $[A]_I$'s with the set-theoretic operations 
modulated by $I:\ [A]_I\lor [B]_I=[A\cup B]_I,\ [A]_I\land [B]_I=[A\cap B]_I,
-[A]_I=[\kappa - A]_I$, ${\bf 0}_I=[\emptyset ]_I$, and
 ${\bf 1}_I=[\kappa ]_I$. The
 subscript $I$ will be suppressed when clear from the context. A Boolean 
algebra is {\it complete iff\/} least upper bounds exist for any collection of
 its members. ${\cal A}$ is an {\it antichain} with respect to $I$ {\it iff} 
${\cal A}\subseteq I^{+}$ yet whenever $A,B\in {\cal A}$ are distinct, $A\cap B
\in I$. Thus, $I$ is $\lambda$-saturated {\it iff} every antichain with respect
 to $I$ has cardinality $<\lambda$. Also, ${\cal A}$ is a maximal antichain 
{\it iff\/} 
$\bigvee\{ [A]\ |\ A\in {\cal A}\}={\bf 1}$ yet whenever $A,B\in {\cal A}$ are
 distinct, $[A]\land [B]={\bf 0}$. 

The first significant results bearing on Sikorski's question were derived by
 Solovay [So], who established the consistency, relative to the existence of a
 measurable cardinal, of the existence of saturated ideals over accessible 
cardinals. In passing, he in effect
 noted the following partial converse to the 
Smith-Tarski result:
\bigskip
{\bf 1.1 Lemma.} {\it Suppose that $I$ is an ideal
 over $\kappa$ such that ${\cal P}
(\kappa )/I$ is complete. For any $\lambda$, if $I$ is not $\lambda$-saturated,
 then $2^\lambda\leq 2^\kappa$. In particular,

(a) $I$ is $2^\kappa$-saturated.

(b) If $2^\kappa <2^{\kappa^{+}}$, then $I$ is $\kappa^{+}$-saturated.}
\medskip
{\bf Proof.} Let $\langle A_\alpha\ |\ \alpha <\lambda\rangle$ enumerate 
(without repetitions) an antichain with respect to $I$. For any $X\subseteq
\lambda$, let $a_X = \bigvee \{ [A_\alpha ]\ |\ \alpha\in 
X\}$. Then $X\not= Y$ implies that $a_X\not= a_Y$. Hence, $2^\lambda\leq\ |\ 
{\cal P}(\kappa )/I|\ \leq 2^\kappa$. \hfill $\dashv$
\bigskip
Kunen established that if there is a $\kappa^{+}$-saturated ideal over 
$\kappa$, then $\kappa$ is measurable in an inner model. In particular, as 
Solovay noted, 1.1(b) implies that if there is an ideal $I$ over $\kappa$ such
 that ${\cal P}(\kappa )/I$ is complete, then $V\not= L$. Kunen asked in the
 early 1970's whether completeness has strong consistency strength, and Theorem
 A confirms this, at least if $\kappa\geq\omega_3$.

With our ultimate goal
 the proof of Theorem A, we now fix an ideal $I$ over $\kappa$ 
such that ${\cal P}(\kappa )/I$ is complete for the rest of this section. We
 use the well-known strategem of considering ${\cal P}(\kappa )/I - \{ 0\}$ 
as a notion of forcing with $[X]\leq [Y]$ {\it iff} $X - Y\in I$, and we 
denote the corresponding forcing relation simply by $\vd$. Note that if
 $\langle A_\alpha\ |\ \alpha <\gamma\rangle$ enumerates without repetitions a
 maximal antichain with respect to $I$, then it corresponds to a name 
$\dot\tau$ for an ordinal specified by: $[A_\alpha ]\ \vd\ \dot\tau =\beta$
{\it iff\/} $\alpha =\beta$.

The following lemmata derive consequences of completeness using maximal 
antichains.
\bigskip
{\bf 1.2 Lemma.} {\it Suppose that $\gamma \leq\kappa$ and 
$\dot\tau$ is a name such
 that $\vd\ \dot\tau \in\gamma$. Then there is a partition $\{ B_\xi\ |\ \xi
 <\gamma\}$ of $\kappa$ such that: if $B_\xi\in I^{+}$, then 
$[B_\xi ]\ \vd\ \dot\tau =\xi$.
}

{\bf Proof.} Let ${\cal A}$ be a maximal antichain with respect to $I$ such
 that for any $A\in {\cal A}$, there is a $\xi_A <\gamma$ such that $[A]
\ \vd\ \dot\tau =\xi_A$. By completeness, let $A_\xi\subseteq\kappa$ for $\xi
 <\gamma$ be such that $[A_\xi ]
=\bigvee\{ [A]\ |\ \xi_A=\xi\}$. Then set $B_\xi =
A_\xi - \bigcup_{\zeta <\xi}A_\zeta$ for $\xi <\gamma$, so that $[B_\xi ]=
[A_\xi ]$ by $\kappa$-completeness. The $B_\xi$'s are as required, once they
 are slightly modified to constitute a partition of all of $\kappa$. 
\hfill $\dashv$
\bigskip 
{\bf 1.3 Lemma.} {\it Forcing with $\vd$ preserves all cardinals $<\kappa$.}

{\bf Proof.} It suffices to show that if $\gamma <\kappa$ is regular, $\delta <
\gamma$, and $\vd\ \dot\tau :\ \delta\rightarrow\gamma$, then $\vd\ 
\exists\eta <\gamma (\dot\tau^{``}\delta \subseteq\eta )$.
For each $\beta <\delta$, let $\langle B_\xi^\beta\ |\ \xi <\gamma\rangle$ 
satisfy 1.2 for $\dot\tau (\beta )$. For each $\alpha <\kappa$, set 
$$
\eta_\alpha =\sup\{ \xi <\gamma\ |\ \alpha\in B_\xi^\beta\ \hbox{\rm and}\ 
\beta <\delta\}\;,
$$ 
so that $\eta_\alpha <\gamma$ by the regularity of $\gamma$.
 Next, set $E_\eta =\{ \alpha <\kappa\ |\ \eta_\alpha =\eta\}$ for $\eta <
\gamma$. Then $\bigcup_{\eta <\gamma}E_\eta =\kappa$ is a partition. 
Consequently, for any $X\in I^{+}$ there is an $\eta <\gamma$ such that $E_\eta
 \cap X\in I^{+}$ by $\kappa$-completeness and $[E_\eta\cap X]
\ \vd\ \dot\tau^{``} \delta\subseteq\eta$.\hfill $\dashv$

The proof of the following proposition
 is similar; it will not be needed in the rest
 of the paper. The only early result about complete quotient Boolean algebras
 other than the Smith-Tarski result, it appeared in terms of distributivity
 properties in Pierce [P], which also contained a similar formulation of the
 easy forcing fact that a notion of forcing adjoins a new function: $\lambda
\rightarrow\lambda$ exactly when it adjoins a new subset of $\lambda$.
\bigskip
{\bf 1.4 Proposition.}
 {\it Suppose that $2^\nu <\kappa$. Then forcing with $\vd$
 does not adjoin any new functions: $\nu\rightarrow 2^\nu$.}

{\bf Proof.} Suppose that $\vd\ \dot\tau :\ \nu\rightarrow 2^\nu$. For each 
$\beta <\nu$, let $\langle B_\xi^\beta\ |\ \xi <2^\nu\rangle$ be as in 1.2 with
 its $\gamma =2^\nu$ and
 its $\dot\tau =\dot\tau (\beta )$. For each $f:\nu\rightarrow
 2^\nu$, let $E_f=\bigcap_{\beta <\nu}B_{f(\beta )}^\beta$. Then $\bigcup_f 
E_f=\kappa$ is a partition. Consequently, for any $X\in I^{+}$, there is an $f$
 such that $E_f\cap X\in I^{+}$ by $\kappa$-completeness, and $[E_f\cap X]\ 
\vd\ \dot\tau =\check f$.\hfill $\dashv$
\bigskip
The connection with inner models of measurability is made through the 
well-known concept of precipitous ideal, due to Jech and Prikry. For an ideal
 $J$ over $\lambda$, if $G$ is generic
 over $V$ for the corresponding notion of forcing
 ${\cal P}(\lambda )/J - \{ 0_J\}$, then $\{ X\subseteq\lambda\ |\ [X]\in
 G\}$ is an ultrafilter on ${\cal P}(\lambda )\cap V$, and
 for any class $A$ in the sense of $V$,
 the ultrapower of
 $A$ with respect to this ultrafilter using functions in  
${}^\lambda V\cap V$ is called the {\it generic ultrapower of} $A$ by $G$.
 In this situation, if $\dot\tau$ is a name for a function in 
${}^\lambda A\cap V$, we denote by $(\dot\tau )$ the name of the equivalence 
class of 
$\dot\tau$ in the ultrapower. $J$ is {\it precipitous iff\/} for any such $G$,
 the generic ultrapower of $V$ by $G$ is well-founded.
 As Jech and Prikry showed, {\it if there is a precipitous ideal $J$ over 
$\lambda$, then $\lambda$ is measurable in an inner model}.

Before we turn to the proof of Theorem A, we establish a partial 
well-foundedness result about generic ultrapowers. 
Continuing to work with our fixed ideal $I$ over $\kappa$ such that ${\cal P}
(\kappa )/I - \{ 0\}$ is complete and the corresponding notion for forcing
 $\vd$, we first make an observation related to 1.2.
\bigskip
{\bf 1.5 Lemma.} {\it Suppose that $\gamma <\kappa$ and $[X]\ \vd\ \dot\tau
\colon \kappa\rightarrow\gamma\ \land\ \dot\tau\in \check V$. Then there is an
 $f\colon\kappa\rightarrow\gamma$ such that $[X]\ \vd\ (\dot\tau )=(\check f)$,
 i.e.~$[X]$ forces that $\dot\tau$ and $\check f$ belong to the same 
equivalence class in the generic ultrapower.}
\medskip
{\bf Proof.} For any set $x$, let $c_x$ denote the constant function: $\kappa
\rightarrow \{ x\}$. We first note that for any $[Y]\leq [X]$, there is a $[Z]
\leq [Y]$ and a $\xi <\gamma$ such that $[Z]\ \vd\ (\dot\tau )=
(\check c_\xi )$. (There is a $[\overline{Y}]\leq [Y]$ and a $g\colon\kappa
\rightarrow\gamma$ such that $[\overline{Y}]\ \vd\ \dot\tau =\check g$. By
 $\kappa$-completeness there is a $\xi <\gamma$ such that $Z=\{
 i\in\overline{Y}\ |\ g(i)=\xi\}\in I^{+}$. Thus, $[Z]\ \vd\ (\dot\tau )=
(\check c_\xi )$ by definition of generic ultrapower.)

Using this, let ${\cal A}\subseteq {\cal P}(X)\ \cap\ I^{+}$ be a maximal 
antichain in ${\cal P}(X)\ \cap\ I^{+}$ such that for each $A\in {\cal A}$, 
there is a $\xi_A <\gamma$ with $[A]\ \vd\ (\dot\tau )=(\check c_{\xi_A})$. By
 completeness, let $A_\xi\subseteq X$ for $\xi <\gamma$ be such that $[A_\xi ]
=\bigvee \{ [A]\ |\ \xi_A=\xi\}$. Then as in the proof of 1.2, let $B_\xi
\subseteq X$ for $\xi <\gamma$ be such that $[B_\xi ]=[A_\xi ]$ and $\{ B_\xi\ 
|\ \xi <\gamma\}$ is a partition of $X$. Finally, define $f\colon\kappa
\rightarrow\gamma$ by $f(i)=\xi$\ {\it iff}\ $i\in B_\xi$ (and $f(i)=0$ if
 $i\not\in X$). Then $[X]\ \vd\ (\dot\tau )=(\check f)$ by definition of 
generic ultrapower. \hfill $\dashv$
\bigskip
This leads to:
\medskip
{\bf 1.6 Proposition.}
 {\it If $\kappa$ is a successor cardinal, then $\vd$ ``the generic
 ultrapower of $\kappa$ is well-founded.''}
\medskip
{\bf Remark.} It follows from the conclusion that e.g.~$\k$ is 
inaccessible in $L$.
\medskip
{\bf Proof.} Suppose that $\kappa =\mu^{+}$, and assume to the contrary that 
for some
 $X\in I^{+},\ [X]\ \vd\ ``\langle (\dot\tau_n)\ |\ n\in\omega\rangle$ is
 an infinite descending sequence with $\dot\tau_n\in {}^\kappa\kappa\ \cap\ 
\check V$ for $n\in\omega$''. We can assume that $[X]\ \vd\ (\dot\tau_0)=
(\check g_0)$ for some $g\in {}^\kappa\kappa$, and that $[X]\ \vd\ \forall n\in
\omega\forall\xi <\kappa (\dot\tau_n(\xi )\leq \check g_0(\xi ))$.

For each $\alpha <\kappa$, let $e_\alpha\colon \alpha +1\rightarrow\mu$ be an
 injection. Then let $\dot\sigma_n$ for $0<n<\omega$ be names such that $[X]\ 
\vd\ \forall\xi <\kappa (\dot\sigma_n(\xi )=\check e_{g_0(\xi )}(\dot\tau_n
(\xi )))$. Thus, $[X]\ \vd\ \dot\sigma_n\colon\kappa\rightarrow\mu$. By 1.5,
 let $f_n\colon\kappa\rightarrow\mu$ for $n\in\omega$ be such that $[X]\ \vd\ 
(\dot\sigma_n)=(\check f_n)$. Next, define $g_n\colon\kappa\rightarrow\kappa$
 for $0<n<\omega$ by: $g_n(i)$ is the unique $\beta$ such that $f_n(i)=e_{g_0
(\xi )}(\beta )$ if $f_n(i)$ is in the range of $e_{g_0(\xi )}$, and $=0$
 otherwise. Then for each $n\in\omega,\ [X]\ \vd\ (\dot\tau_n)=(\check g_n)$.
 Finally, if $T_n=\{ i\in\kappa\ |\ g_{n+1}(i)<g_n(i)\}$ for $n\in\omega$, then
 $[X]\leq [T_n]$ by our assumption on $[X]$ and the definition of generic 
ultrapower. Hence, $X\ \cap\ \bigcap_nT_n$ is not empty, but any $i$ in this
 set gives rise to an infinite descending set of ordinals $g_0(i)>g_1(i)>g_2(i)
\ldots$, which is a contradiction.\hfill $\dashv$
\bigskip
Turning to the proof of Theorem A,
 we need another ingredient. The Dodd-Jensen Covering Theorem for their inner
 model $K$ asserts that if there is no inner model with a measurable 
cardinal, then for any uncountable set $x$ of ordinals there is a set $y
\supseteq x$ such that $y\in K$ and $|y|=|x|$. The definable class $K$ is 
extensionally preserved in all forcing extensions.

{\bf Proof of Theorem A.} We take $\kappa\geq\omega_3$. 
Assume to the contrary that there is no inner model
 with a measurable cardinal. Suppose briefly that $G$ is generic for ${\cal P}
(\kappa )/I - \{ 0\}$. Then $K^V=K^{V[G]}$, and (since forcing does not 
create inner models of measurability) the conclusion of the Covering Theorem 
holds in $V[G]$. Since $\kappa\geq\omega_3$ and $\omega_1$ and $\omega_2$ are
 preserved by the forcing by 1.3, this implies in particular: whenever $x$ is
 a size $\omega_1$ set of ordinals in $V[G]$, there is a $y\in V$ of size 
$\omega_1$ such that $y\supseteq x$. (Note that the preservation of $\omega_2$
 is needed here.) Recalling a previous remark about maximal antichains 
corresponding to names for ordinals, this in turn translates in a 
straightforward manner to the following assertion in the ground model:
\medskip
{\parindent=30pt \narrower \noindent \llap {\hbox to \parindent 
{$(\ast )$\hfil}}Whenever $\langle {\cal A}_\alpha\ |\ \alpha <
\omega_1\rangle$ is a sequence
 of maximal antichains with respect to $I$, for any $X\in I^{+}$ there 
is a $Y\in {\cal P}(X)\cap I^{+}$ such that for any $\alpha <\omega_1$,
 $|\{ A\in {\cal A}_\alpha\ |\ A\cap Y\in I^{+}\}|\leq\omega_1$.\par}
\medskip
We now derive a contradiction by using this to show that $I$ must
 be precipitous. Suppose then that $[X]\ \vd\ ``\langle (\dot\tau_n)
\ |\ n\in\omega
\rangle$ is an infinite descending sequence in the
 generic ultrapower''.
 For each $n\in\omega$, let ${\cal A}_n$ be
 a maximal antichain with respect to $I$ such that whenever $A\in {\cal A}_n$
 there is an $f_n^A:\kappa\rightarrow V$ such that $[A]\ \vd\ \dot\tau_n=
\check f_n^A$. Let $Y\in {\cal P}(X)\cap I^{+}$ satisfy $(\ast )$ for $\langle
 {\cal A}_n\ |\ n\in\omega\rangle$ and for each $n\in\omega$ let $\langle 
A_\xi^n\ |\ \xi <\omega_1\rangle$ enumerate $\{ A\in {\cal A}_n\ |\ A\cap Y
\in I^{+}\}$. Set $B_\xi^n=(A_\xi^n - \bigcup_{\zeta <\xi}A_\zeta^n)\cap Y$
 and $B_n=\bigcup_{\xi <\omega_1}B_\xi^n$, and define $f_n$ on $B_n$ by: $f_n=
\bigcup_{\xi <\omega_1}f_n^{A_\xi^n}|B_\xi^n$. It is easily seen that for $n\in
\omega,\ [B_n]=[Y]$ by maximality of ${\cal A}_n$ and $\kappa$-completeness.
 Hence $\bigcap_nB_n\not=\emptyset$, again by $\kappa$-completeness. But for
 any $i$ in this set, $f_0(i)>f_1(i )>\ldots$, which is a 
contradiction.\hfill $\dashv$ 
\bigskip
Properties of the sort $(\ast )$ were first investigated 
in Baumgartner-Taylor [BT], then in Foreman-Magidor-Shelah [FMS], by Woodin, 
and
extensively in Gitik-Shelah [GS]. A positive answer to the following
 question would strengthen 1.6 and eliminate the condition $\kappa\geq\omega_3$
 from Theorem A.
\medskip
{\bf 1.7 Question}. If $\kappa$ is a successor cardinal and $I$ is an ideal 
over $\kappa$ such that ${\cal P}(\kappa )/I$ is complete, is $I$ precipitous?
\vskip 0.3in
\centerline{$\S$2 Separating Completeness from Saturation}
\medskip
We next turn to the proof of Theorem B.  We shall build on the proof of
 the following result of Shelah [Sh4:~XVI, Theorem 2.4],
 which we first describe. 
\bigskip
{\bf 2.1 Theorem.} {\it Suppose that $\kappa$ satisfies $(\ast )_{ab}^a$ (see
 [Sh4:~XVI, 2.2]), e.g.~is a Woodin cardinal, and $S$ is
 a stationary, co-stationary subset of\/ $\omega_1$. Then there is a forcing 
extension with no new reals satisfying: $\kappa =\omega_2$ and 
${\rm NS}|\tilde S$ is an ideal over $\omega_1$ which is $\omega_2$-saturated.}

This result weakened the large cardinal hypotheses of previous results drawing
 the same conclusion, and was an outgrowth of Foreman-Magidor-Shelah [FMS]. It
 was established using concepts and techniques of the second author that we
 quickly review:

Suppose that $\langle P,\leq\rangle$ is a poset for forcing, $\lambda$ is
 regular with ${\cal P}(P)\subseteq H_\lambda$, and $N$ is countable with 
$\langle N,\in\rangle\prec\langle H_\lambda,\in\rangle$. Then $q\in P$ is 
$\langle N,P\rangle$-{\it generic\ iff}\ for any $P$-name $\dot\tau\in N$ for
 an ordinal, $q\ \vd\  \dot\tau\in\check N$. Refining this,
 $q\in P$ is 
$\langle N,P\rangle$-{\it semigeneric\ iff}
 for any $P$-name $\dot\tau\in N$ for
 a {\it countable} ordinal, $q\ 
\vd\ \dot\tau\in \check N$. $\langle p_n\ |\ n\in
\omega\rangle$ is a $P$-{\it generic sequence for $N$} {\it iff\/} $p_{n+1}\leq
 p_n\in N$ for each $n$ and whenever $D$ is a dense set for $P$ with $D\in N$,
 there is an $n$ such that $p_n\in D$.

$\langle P,\leq\rangle$ is {\it semiproper iff\/} for any regular 
$\lambda$ such that ${\cal P}(P)\subseteq H_\lambda$, there is a closed 
unbounded subset of $[H_\lambda ]^{<\omega_1}$ consisting
 of $N$ such that $\langle N,\in\rangle\prec
 \langle H_\lambda,\in\rangle$ satisfying: for any $p\in N$, there is a $q\leq
 p$ such that $q$ is $\langle N,P\rangle$-semigeneric. For $S$ a stationary 
subset of $\omega_1,\ \langle P,\leq\rangle$ is $S$-{\it closed iff\/} for any
 regular $\lambda$ such that ${\cal P}(P)\subseteq H_\lambda$, countable $N$
 with $\langle N,\in\rangle\prec\langle H_\lambda,\in\rangle$ and $N\cap 
\omega_1\in S$, and $P$-generic
 sequence $\langle p_n\ |\ n\in\omega\rangle$ for
 $N$, there is a $q\leq p_n$ for every $n\in\omega$. Semiproper is $\{ 
\aleph_1\}$-semiproper and $S$-closed is $\{ S\}$-complete in the sense of 
Shelah [Sh2].  $\omega_1$-closure readily implies semiproperness.
The salient features of these concepts are that if a poset 
is semiproper, then forcing with it {\it preserves stationary subsets
 of} $\omega_1$ (i.e. any stationary subset of $\omega_1$ in the ground model
 remains stationary in the extension), and if it is $S$-closed, then it adjoins
 no new countable sequences of ordinals.

To affirm notation,
$P$ is the {\it countable support iteration} of $\langle P_\alpha,\dot Q_\alpha
\ |\ \alpha <\gamma\rangle$ {\it iff\/} setting $P_\gamma =P$, we have: 
$P_0=\{ \emptyset\}$; for $\alpha <\gamma$,\ $\vd_{P_\alpha} ``\dot Q_\alpha$
 is a poset''
 and $P_{\alpha +1}=P_\alpha\ast\dot Q_\alpha$; and for limit $\alpha\leq
\gamma$, $P_\alpha$ is the direct limit of $\langle P_\beta\ |\ \beta <\alpha
\rangle$ in case cf$(\alpha )>\omega$, and the inverse limit otherwise. 
Proceeding recursively we can take $P_\alpha$ to consist of functions
 $p$ (the conditions) with domain $\alpha$ so that for each $\beta <\alpha,\ 
p(\beta )$ is a $P_\beta$-name and $\vd_{P_\beta}\ p(\beta )\in\dot Q_\beta$,
 and supposing that $\vd_{P_\beta}\ ``\dot 1_{Q_\beta}$ is the maximum 
element of $\dot Q_\beta$'',
$$
{\rm supt}(p)=\{ \beta <\alpha\ |\ \vd_{P_\beta}
\ p(\beta )\not= \dot 1_{Q_\beta}\}
$$
is countable, with corresponding partial order on $P_\alpha$ given by: $p\leq
 q$\ {\it iff}\/ $\forall\beta <\alpha (p|\beta\ \vd_{P_\beta}\ p(\beta )\leq 
q(\beta ))$.

Throughout the paper we rely on the following convention: For a notion of
forcing $P$, $\dot G_P$ denotes the canonical $P$-name for its generic
object, and if $P = P_\alpha$ in some contextually clear indexing, 
$\dot G_\alpha$ is written for $\dot G_{P_\alpha}$.
\bigskip
{\bf 2.2 Lemma.} {\it Suppose that $P$ is the countable support iteration of 
$\langle P_\alpha,\dot Q_\alpha\ |\ \alpha <\gamma\rangle$, where for each 
$\alpha <\gamma$, $\vd_{P_\alpha}$\ ``$\dot Q_\alpha$ is semiproper
and $S$-closed''. 
 Then $P$ is semiproper and $S$-closed.}
\bigskip
This is a special case of more general iteration lemmas. The appropriate mode
 of iteration for semiproperness is revised countable support (RCS) iteration,
 but $S$-closure at each stage implies that there are no new countable 
sequences of ordinals, and so RCS iteration reduces to countable support 
iteration.

In [FMS], 2.1 is established for $\kappa$ supercompact instead of,
 e.g.~Woodin, by first establishing the consistency with CH of a Martin's 
Axiom for $S$-closed notions of
 forcing that preserve stationary subsets of $\omega_1$ and meeting $\omega_1$
 dense sets. Then it suffices to argue with a notion of forcing which in the
 formulation of Shelah [Sh3] is as follows:

Suppose that ${\cal A}$ consists of stationary subsets of $\omega_1$ with $S
\in {\cal A}$. Then
$$
{\rm Seal}_S({\cal A})
$$
consists of countable sequences $\langle N_\xi\ |\ \xi\leq\gamma\rangle$ such
 that:

(i) each $N_\xi\in [H_\lambda ]^{<\omega_1}$, where $\lambda =
(2^{\omega_1})^{+}$, and $S\in N_0$.

(ii) $\langle N_\xi\ |\ \xi\leq\gamma\rangle$ is increasing and continuous.

(iii) for each $\xi\leq\gamma$, $N_\xi\cap\omega_1$ is an ordinal in $\bigcup
 \{ A\ |\ A\in {\cal A}\cap N_\xi\}$.

\noindent ${\rm Seal}_S({\cal A})$ is ordered by: $p\leq q$ {\it iff} $q$ is 
an initial segment of $p$.

${\rm Seal}_S({\cal A})$ is clearly $S$-closed, although it may not be 
semiproper,
 and forcing with it provides an enumeration of ${\cal A}$ in ordertype $\leq
\omega_1$ and a closed unbounded subset of the diagonal union of ${\cal A}$
 according to that enumeration. If ${\cal A}$ was a maximal antichain with
 respect to NS, then it can be shown that forcing with ${\rm Seal}_S({\cal A})$ 
preserves stationary subsets of $\omega_1$. ${\cal A}$ is then ``sealed'': it
 remains a maximal antichain in any extension that preserves stationary sets,
 since any stationary set has stationary intersection with the diagonal union
 of ${\cal A}$, and hence with a particular member of ${\cal A}$ by Fodor's
 Lemma. The aforementioned version of Martin's Axiom implies through this means
 that ${\rm NS}|\tilde S$ is $\omega_2$-saturated.

2.1 was established by applying reflection properties 
directly in [Sh4:~XVI]. Let
$$
Q_S=({\textstyle \prod} \{ {\rm Seal}_S({\cal A})\ |\ {\rm Seal}_S({\cal A})\ 
\hbox{\rm is semiproper}\})\ast {\rm col}(\omega_1,2^{\omega_1})\;,
$$
the countable support product of all 
${\rm Seal}_S({\cal A})$'s for ${\cal A}$'s that yield semiproper 
${\rm Seal}_S({\cal A})$, followed by the usual collapse of
 $2^{\omega_1}$ to $\omega_1$ using countable approximations. Clearly $Q_S$ 
is   $S$-closed; it is also semiproper (see [Sh4:~XIII, 2.8(3)] or 
[Sh3:~2.8(3), p.~361]). Let $P$ be the countable support iteration of
 $\langle P_\alpha,\dot Q_\alpha\ |\ \alpha <\kappa\rangle$ where 
$\dot Q_\alpha$ is a $P_\alpha$-name for $Q_S$ in the sense of $V^{P_\alpha}$.
 Assuming that $\kappa$ satisfies the
 large cardinal hypothesis of 2.1, the second author showed that any forcing
 extension via $P$ satisfies its conclusion.

{\bf Proof of Theorem B.} 
 We interlace into the above described proof of 2.1 
natural notions of forcing for introducing suprema into ${\cal P}(\omega_1)/
{\rm NS}|\tilde S$: $\langle {\cal A},{\cal B}\rangle$ is an {\it appropriate
 pair} if ${\cal A}$ and ${\cal B}$ consist of  subsets of $\omega_1$
 such that: if $A\in {\cal A}$ and $B\in {\cal B}$, then $A\cap B\cap\tilde S$
 is non-stationary. For such $\langle {\cal A},{\cal B}\rangle$, let 
${\rm Sup}_S({\cal A},{\cal B})$ consist of triples $\langle w,c,d\rangle$ 
such that:

(i) $w$ is a countable subset of ${\cal A}\cup {\cal B}$.

(ii) $c$ is a countable, closed set of countable ordinals (so $\cup 
c\in c$).

(iii) $d: c\rightarrow 2$.

Order ${\rm Sup}_S({\cal A},{\cal B})$ by: $\langle w,c,d\rangle\leq\langle 
\overline{w},\overline{c},\overline{d}\rangle$ {\it iff\/} 
$\overline{w}\subseteq w,\ \overline{c}$ is an initial segment of 
$c,\ \overline{d}\subseteq d$, and:
\medskip
\centerline{
if $\alpha\in (c\sim\overline{c})\cap\tilde S$ and $d(\alpha )=0$, then $\alpha
\notin A$ for any $A\in\overline{w}\cap {\cal A}$.}
\centerline{
if $\alpha\in (c\sim\overline{c})\cap\tilde S$ and $d(\alpha )=1$, then $\alpha
\notin B$ for any $B\in \overline{w}\cap {\cal B}$.}
\medskip
Suppose that $G$ is ${\rm Sup}_S({\cal A},{\cal B})$-generic and set 
$c_G=\bigcup \{
 c\ |\ \langle w,c,d\rangle\in G\}$ and $d_G=\bigcup \{ d\ | \ \langle w,c,d
\rangle\in G\}$. Then $c_G$ is a closed unbounded subset of $\omega_1$, 
$d_G^{-1}(\{ 0\})\cap A\cap\tilde S$ is
 countable for every $A\in {\cal A}$, and
 $d_G^{-1}(\{ 1\})\cap B\cap\tilde S$ is countable for every $B\in {\cal B}$.
 In particular, if ${\cal A}\cup {\cal B}$ were a maximal antichain with 
respect to ${\rm NS}|\tilde S$, then $[d_G^{-1}(\{ 1\})]$ would serve as an 
upper bound of ${\cal A}$ and $[d_G^{-1}(\{ 0\})]$ of ${\cal B}$.
These will be {\it least} upper bounds in the extension if $\cal A \cup
\cal B$ continues to be maximal there, and this is the only situation
that will be germane to the overall argument. Of course, for all this
to make sense in the extension we must ascertain that stationary subsets of 
$\omega_1$ are preserved:
\bigskip
{\bf 2.3 Lemma.}

{\it (a) ${\rm Sup}_S({\cal A},{\cal B})$ is $\omega_1$-closed (and hence
semiproper and $S$-closed).

(b) Assuming {\rm CH}, ${\rm Sup}_S({\cal A},{\cal B})$ is 
$\omega_1$-linked with least upper bounds, i.e. it is the union of $\omega_1$ sets, each consisting of 
pairwise compatible elements with least upper bounds.}

{\bf Proof.} (a) Set $P={\rm Sup}_S({\cal A},{\cal B})$. Suppose that $\lambda$ is
 regular with ${\cal P}(P)\subseteq H_\lambda,\ N$ is countable with
 $\langle N,\in\rangle\prec\langle H_\lambda,\in\rangle$ and $\langle p_n\ |\ 
n\in\omega\rangle$ is a $P$-generic sequence for $N$. We must find a $q\in P$ 
such that $q\leq p_n$ for every $n\in\omega$.

For $n\in\omega$ and $p_n=\langle w_n,c_n,d_n\rangle$, since $w_n$ is 
countable, there is a $C_n\in N$
be a closed unbounded subset of $\omega_1$ such that for any $A\in
w_n \cap {\cal A}$ and $B\in w_n \cap {\cal B}$, $C_n\cap A\cap
B\cap\tilde S=\emptyset$. Set
$w=\bigcup w_n$, $c=\bigcup c_n$, $d=\bigcup d_n$. Then a simple
genericity argument implies that $\cup c\in \cap C_n$ so that for no
$A\in w\cap {\cal A}$ and $B\in w\cap {\cal B}$ does $\cup c\in A\cap
B$, and so we can find an $i<2$ such that $\langle w,c\cup\{ \cup
c\},d\cup\{ \langle \cup c,i\rangle\}\rangle$ is in ${\rm Sup}_S({\cal A},
{\cal B})$ and of course is $\leq p_n$ for every $n\in\omega$.

(b) Note that $\langle w_0,c,d\rangle$ and $\langle w_1,c,d\rangle$ are 
compatible with least upper bound $\langle w_0\cup w_1,c,d\rangle$. With CH,
 there are $\omega_1$ such pairs.\hfill $\dashv$
\bigskip 
Let
$${\rm Sup}_S$$
be the countable support product of ${\rm Sup}_S({\cal A},{\cal B})$ for all
appropriate pairs $\langle {\cal A},{\cal B}\rangle$. A countable support
product of $\omega_1$-closed forcings is readily seen to be 
$\omega_1$-closed, and
$\omega_1$-closed forcings are $S$-closed and semiproper. Hence by 2.3(a),
{\it ${\rm Sup}_S$ is $S$-closed and semiproper.}

 An approach to the proof of Theorem B would bec to carry out the countable 
support iteration of ${\rm Sup_S}$ through $\kappa$ stages. Just assuming
cf$(\kappa)>\omega_1$, we would then get the consistency of
$2^{\omega_1}=\kappa$ and every appropriate pair $\langle {\cal A}, 
{\cal B}\rangle$ with $| {\un A} \cup {\un B}| < \k$ can be separated, a 
consequence in fact of a 
generalized Martin's Axiom in Baumgartner [B] or in Shelah [Sh1].  It
is to ensure $\kappa$-saturation, a necessary condition for full
completeness by 1.1(a), that we build on the proof of 2.1.

Let $P_\kappa$ be the countable support iteration of $\langle P_\alpha, 
\dot Q_\alpha\ |\ \alpha <\kappa\rangle$, where:

(i) For odd $\alpha <\kappa,\ \dot Q_\alpha$ is a $P_\alpha$-name for $Q_S$ in
 the sense of $V^{P_\alpha}$. (Here, $Q_S$ is as in the above outline of the
 proof of 2.1.)

(ii) For even $\alpha <\kappa,\ \dot Q_\alpha$ is a $P_\alpha$-name for the 
countable support product of ${\rm Sup}_S({\cal A},{\cal B})$'s for all 
appropriate pairs $\langle {\cal A},{\cal B}\rangle$ satisfying 
$\langle {\cal A},
{\cal B}\rangle\in V[\dot G_\alpha|\{ 2\gamma\ |\ 2\gamma <\alpha\}]$. (For 
$Z\subseteq\alpha$, $\dot G_\alpha|Z$ is the $P_\alpha$-name for 
$\{ p\in\dot G_\alpha\ |\ {\rm supt}(p)\subseteq  Z\}$. Note that
$\langle {\un A}, {\un B} \rangle $ is to belong to the smaller model, but 
in the definition of ``appropriate pair'' the non-stationariness in 
``$(\forall A \in {\un A}) (\forall B \in {\un B}) \; A \cap B \cap \tilde S$
is not stationary'' is to be in the sense of $V [G_\alpha] $!) 

Next, for $\alpha\leq\kappa$ set, by induction on $\alpha$:
$$
P_\alpha^\prime =\{ p\in P_\alpha\ |\ {\rm supt}(p)\subseteq\{ 2\gamma\ |\ 
2\gamma <\alpha\}\ \land\ \forall \beta \in {\rm supt}(p)(p(\beta)\ \hbox{\rm is a}\ P^\prime_{\beta}{\rm -name})\}
$$
with the inherited order. ($p(\beta )$ being a $P^\prime_\beta$-name, it
only depends on $\dot G_\beta \cap P^\prime_\beta$.)  We show 
that any forcing extension via
$P_\kappa^\prime$ satisfies the conclusion of Theorem B. This follows
from the following technical lemma, all of whose parts are established
by simultaneous induction; for its (b), note that $P^\prime_\beta$-names
being $P_\beta$-names is justified by an inductive appeal to
(a), and ${\rm Sup}_S$ was defined a few paragraphs ago in the outline
of the proof of 2.1.
\bigskip
{\bf 2.4 Lemma.} {\it For each $\alpha\leq\kappa$:

(a) $P_\alpha^\prime \lessdot   P_\alpha$, i.e. every maximal antichain of 
$P_\alpha^\prime$ is a maximal antichain of $P_\alpha$.

(b) $P_\alpha^* =\{ p\in P_\alpha\ |\ \forall {\it even}\
\beta\in {\rm supt}(p)(p(\beta )\ \hbox{\it is a}\
P_\beta^\prime\hbox{\it -name for a condition}$ in ${\rm Sup}_S$  in the 
sense of\/ $V^{P_\beta^\prime}\}$ \hbox{\it is dense in} $P_\alpha$. 

(c) For any $p,q\in P_\alpha^*$ such that $p|\{ 2\gamma\ | \ 2\gamma
<\alpha\}$ and $q |\{ 2\gamma\ |\ 2\gamma <\alpha\}$ are compatible members
of $P^\prime_\alpha$, there is an automorphism
$F_{pq}^\alpha$ of $P_\alpha$ such that: $F_{pq}^\alpha (p)$ is
compatible with $q$; $F_{pq}^\alpha$ is the identity on
$P_\alpha^\prime$; and inductively for any $\beta <\alpha$,
$$F_{p|\beta, q|\beta}^\beta =F_{pq}^\alpha |P_\beta .$$

(d) If $p \in P^*_\alpha$, then $p |  \{ 2\gamma\ |\ 2\gamma < \alpha\}
\in P^\prime_\alpha$. 
}
\bigskip
Once this lemma is established, the proof of Theorem B can be
completed as follows: Suppose that $G^\prime$ is $P_\kappa^\prime$-generic,
and by 2.4(a) let $G$ be $P_\kappa$-generic such that $G\cap P_\kappa^\prime
= G^\prime$. Note first that for any $X\subseteq\omega_1$ with $X
\in V[G^\prime ]$,
$$
V[G^\prime ]\models  X\ \hbox{\rm is stationary}\ \ {\it iff}\ \ V[G]\models
X\ \hbox{\rm is stationary}\;.
$$
(If $X$ were non-stationary in $V[G]$, then it would be non-stationary in 
$V[G\cap P_\alpha ]$ for some even $\alpha <\kappa$. But then, $X=A\cap B\cap
\tilde S$ for some $A\in {\cal A}$ and $B\in {\cal B}$ with $\langle {\cal A},
{\cal B}\rangle$ an appropriate pair in $V[G\cap P_\alpha^\prime ]$, so that
 ${\rm Sup}_S({\cal A},{\cal B})$ at that stage would have adjoined a closed 
unbounded subset of $\omega_1$ confirming that $X$ is non-stationary in 
$V[G^\prime ]$.) The proof of 2.1 still works to show that in $V[G],\ \kappa =
\omega_2$ and ${\rm NS}|\tilde S$ is $\kappa$-saturated. It thus follows that 
${\rm NS}|\tilde S$ is also $\kappa$-saturated in $V[G^\prime ]$.

We can conclude that for any maximal antichain ${\cal A}\cup {\cal B}$ of 
${\rm NS}|\tilde S$ in the sense of $V[G^\prime ]$ with $\cal A \cap \cal B 
= \emptyset$, $\langle {\cal A},{\cal B}\rangle\in V[G^\prime
\cap P_\alpha^\prime ]$ for some even $\alpha <\kappa$. But then, the forcing
with ${\rm Sup}_S ({\cal A},{\cal B})$ would have adjoined a set $E \subseteq 
\omega_1$ such that: if $A \in \cal A$, then $A - (E \cap \tilde S)$ is not 
stationary; and if $B \in \cal B$, then  $B \cap E \cap \tilde S$ is not 
stationary. $[E \cap \tilde S]$ is thus the supremum for ${\cal A}$ in 
$V[G^\prime\cap P_{\alpha +1}^\prime ]$ and 
hence in $V[G^\prime ]$. This suffices to show that in $V[G^\prime ]$, 
${\cal P}(\omega_1)/{\rm NS}|\tilde S$ is a complete Boolean algebra.

Finally, $P_\kappa^\prime$ is $S$-closed, so that forcing with it
adjoins no new countable sequences of ordinals. Consequently, 
\medskip
\centerline{$\{ p\in P_\kappa^\prime\ |$ if $\beta\in {\rm supt}(p)$ and 
$\langle \dot w,\dot c,\dot d\rangle$ is a component of $p(\beta )$, then $\dot
c=\check c$ and $\dot d=\check d$ for some $c,d\in V\}$}
\medskip\noindent
is dense in $P_\kappa^\prime$. Using CH and 2.3(b), a simple $\Delta$-system
argument using this dense set then shows that $P_\kappa^\prime$ has
the $\omega_2$-c.c. and hence preserves all cardinals, and it is
simple to see that it renders $2^{\omega_1}=\kappa$.

{\bf Proof of 2.4}. Assuming that $\alpha\leq\kappa$ and all four
parts hold below $\alpha$, we verify that they all hold at $\alpha$.

To first verify that (b) holds at $\alpha$, suppose that $p\in
P_\alpha$ is arbitrary. Let $\lambda$ be regular and sufficiently
large, and $N$ countable such that $\langle N,\in\rangle\prec\langle
H_\lambda,\in\rangle$, $p\in N$, and $N\cap\omega_1\in S$. Let
$\langle p_n\ |\ n\in\omega\rangle$ be a $P_\alpha$-generic sequence
for $N$ with $p_0=p$. Incorporating the proof of 2.3(a) into the
iteration lemma for $S$-closed notions of forcing, there is a {\it
least} upper bound $q\in P_\alpha$ for the $p_n$'s specified as
follows:

supt$(q)=\bigcup_n {\rm supt}(p_n)$.  For even $\beta\in{\rm
supt}(q)$, $q|\beta$ forces that for component $\langle \dot w,\dot
c,\dot d\rangle$ of $q(\beta )$ in some relevant ${\rm Sup}_S(\dot {\cal
A},\dot {\cal B})$ with corresponding $\langle \dot w_n,\dot c_n,\dot
d_n\rangle$ in $p_n$ for $n\in
\omega$ sufficiently large, $\dot w$ is the union of the $\dot w_n$'s, 
$\dot c$ is the union of the $\dot c_n$'s together with its limit
point at the top, and $\dot d$ is the union of the $\dot d_n$'s
together with an arbitrary value for that top limit point. By
$P_\alpha$-genericity of $\langle p_n\ |\ n\in\omega\rangle$ and
induction it can be assumed that each $\dot w_n$ is a
$P_\beta^\prime$-name and hence that $\dot w$ is a
$P_\beta^\prime$-name, and that $\dot c=\check c$ for some $c\in V$
and $\dot d=\check d$ for some $d\in V$.

Continuing to consider that specific component $\langle \dot w,\dot
c,\dot d\rangle$ of $q(\beta)$, by definition of such conditions there is a
$P_\beta^\prime$-name $\dot w_0$ such that $q|\beta\
\vd\ ``\langle \dot w\cap\dot w_0,\dot w - \dot w_0
\rangle$ is an appropriate pair''. By induction, $P_\beta^\prime 
\lessdot\ P_\beta$, and the homogeneity property 2.4(c) implies that
whenever $r\in P_\beta$, $\dot\tau$ is a $P_\beta^\prime$-name, $\psi$
is a one-free variable formula, and $r\ \vd_{P_\beta}\ \psi (\dot\tau
)$, then $r|\{ 2\gamma\ |\ 2\gamma <\beta\}\
\vd_{P_\beta^\prime}\ \psi (\dot\tau )$. In particular 
$$
(q|\beta )|\{2\gamma\ |\ 2\gamma <\beta\}\ \vd_{P_\beta^\prime}\ ``\langle 
\dot w\cap\dot w_0,\ \dot w - \dot w_0\rangle\ \hbox{\rm is an appropriate
pair''.}$$ 
This confirms that $q\in P_\alpha^*$ to verify 2.4(b) for
$\alpha$ as desired.

At the referee's urging we also elaborate  the rest: 

To establish (d) at $\alpha$, let $p \in P^*_\alpha$, $\beta$ an even 
ordinal in supt$(p)$, and set $p^e = p|\{2\gamma\ |\ 2\gamma < \beta\}$. It
must be shown that $p^e\ \vd_{P_\beta}\ p (\beta) \in \dot 
Q_\alpha$. We know inductively that $P^\prime_\beta \lessdot P_\beta$,
$p^e\in P^\prime_\beta$, and $p^e\ \vd_{P^\prime_\beta}\ p(\beta )\in \dot 
Q_\beta$. The only problem in trying to replace $P^\prime_\beta$ by
$P_\beta$ here is that an appropriate pair mentioned in $p(\beta )$ in the
sense of $V^{P^\prime_\beta}$ may no longer be one in the sense of
$V^{P_\beta}$. 

Assume to the contrary that for some $q\in P_\beta$ with $q\leq p^e$,
$q\ \vd_{P_\beta}\ \dot p (\beta) \notin \dot Q_\beta$. By (b) inductively
it can be assumed that $q \in P^*_\beta$, and by (d) inductively,
$q|\{2\gamma\ |\ 2\gamma <\beta\}\in P^\prime_\beta$. Since
$q|\{2\gamma\ |\ 2\gamma <\beta\}\leq p^e$ in $P^\prime_\beta$, there is
an automorphism $F^\beta_{q,p|b}$ as in (c) inductively such that
$F^\beta_{q,p| \beta} (q)$ is compatible with $p| \beta$. Hence, for some 
$q^+\in P_\beta$ with $q^+ \leq q$, $F^\beta_{q, p| \beta} (q^+) \leq p|\beta$.
Let  $G \subseteq P_\beta$ be $P_\beta$-generic over $V$ with $q^+ \in G$.
Then $G^\prime =F^\beta_{q,p| \beta}``G$ is  $P_\beta$-generic 
over $V$ such that $p| \beta \in G'$.  But $G \cap P^\prime_\beta = 
G^\prime \cap P^\prime_\beta$ as $F^\beta_{q,p| \beta}$ is the identity on
$P^\prime_\beta$, and so
$$ 
\{ \langle{\cal A}, {\cal B}\rangle\ |\ {\cal A} \cup {\cal B} \subseteq
{\cal P} (\omega_1) \cap V (G \cap P^\prime_\beta )\}=\{ \langle{\cal A}, 
{\cal B}\rangle\ |\ {\cal A} \cup {\cal B} \subseteq {\cal P} (\omega_1)\cap 
V [G^\prime \cap P^\prime_\beta] \}\;,
$$ 
and so as $V [ G] =V [G^\prime ])$,
$$\displaylines{
\qquad\{\langle \cala, \calb\rangle \in V[G \cap P^\prime_\beta ]\ |\ \langle 
\cala, \calb\rangle\ \hbox{\rm is appropriate in}\ V [G])\}=\hfill\cr
\hfill
\{\langle\cala,\calb\rangle \in V[G\cap P^\prime_\beta]\ |\ \langle
\cala, \calb\rangle\ \hbox{\rm is appropriate in}\ V[G^\prime]) \}\;.\qquad\cr}
$$
As also all $\omega$-sequences from $V$ of members of $V [G \cap
P^\prime_\beta]$ are in $V [G\cap P^\prime_\beta]$ (as this holds for
$V$ and $V[G]$), clearly $(\dot Q_\beta)^ G= (\dot Q_\beta )^{G^\prime}$.
This contradicts the choice of $q$.

Now clause (a) for $\alpha$ follows: For $p \in P_\alpha$, choose $q$ such that
$q \leq p \in P^*_\alpha$. Setting $q^e = q |\{ 2\gamma\ |\ 2\gamma <
\alpha\}$, $q^e \in P^\prime_\alpha$ by clause (d) and
$q^e\ \vd_{P'_\alpha}\ ``p\in P_\alpha/\dot{G}_{P_\alpha^\prime}$'' because
if $r \leq q^e \in P^\prime_\alpha$ then $r'=r\cup (q|\{\beta < \alpha\
|\ \beta\ {\rm is\ odd}\})$ is in $P^*_\alpha$ (check!) and is below $r$ and 
$q$ hence below $p$. 

We lastly deal with clause (c) for $\alpha$. If $\alpha$ is a limit, it
is immediate: $F^\alpha_{p,q}(r)$ is defined by: $\gamma\in {\rm supt}(
F^\alpha_{p,q}(r))$ {\it iff\/} for some (equivalently, every) $\beta\in 
(\gamma,\alpha)$ we have $\gamma\in {\rm supt}(F^\beta_{p|\beta,q|\beta}(r))$,
and letting $(F^\alpha_{p,q}(r))(\gamma)=(F^\beta_{p|\beta,q|\beta}(r))
(\gamma)$ for some (any) such $\beta$.

If $\alpha=\beta+1$, $\beta$ odd, just note that
$\dot{Q}_\beta$ is definable in ${\rm V}[\dot{G}_{P_\beta}]$ (without
parameters). If $\alpha =\beta+1$, $\beta$ even, if $G\subseteq
P_\beta$ is generic over V, in ${\rm V}[G\cap P^\prime_\beta]$ there
is an automorphism $F_\beta$ of $\dot{Q}_\beta[G]$ (note:
$\dot{Q}_\beta[G]\in V[G\cap P^\prime_\beta]$ -- see the proof of
clause (d)) such that: $F_\beta(p(\beta))$, and $q(\beta)$ are
compatible. (The simplest way to see this is to replace in the iteration
$\dot{Q}_\alpha$ by its completion.)

\vskip 0.3in
\centerline{$\S$3 Collapsing $\kappa$ to $\omega_3$.}
\medskip
In this section we specify the modifications to the proof of Theorem B
necessary to establish Theorem C. We define the components of an
iteration $\langle P_\alpha,\dot R_\alpha\ |\ \alpha <\kappa\rangle$
in three cases instead of two:

(i) For $\alpha\equiv 0$ (mod 3), $\dot R_\alpha$ is defined as $\dot
Q_\alpha$ was before for odd $\alpha$, i.e. it is $Q_S^{V^{P_\alpha}}$.

(ii) For $\alpha\equiv 1$ (mod 3),\ $\dot R_\alpha$ is defined as $\dot
Q_\alpha$ was before for even $\alpha$, but for all appropriate pairs
$\langle {\cal A},{\cal B}\rangle$ satisfying $\langle {\cal A},{\cal
B}\rangle\in V[\dot G_\alpha |\{ \beta <\alpha\ |\ \beta\equiv 1\ ({\rm mod}\
3)\ \vee\ \beta
\equiv 2\ ({\rm mod}\ 3)\}]$ (where $\dot G_\alpha |Z$ is as before).

(iii) For $\alpha\equiv 2$ (mod 3), $\dot R_\alpha$ is col$(\omega_2,
2^{\omega_2})$, the collapse of $2^{\omega_2}$ to $\omega_2$ using
$\omega_1$ size approximations, in the sense of $V[\dot G_\alpha|\{
\beta <\alpha\ |\
\beta\equiv 1\ ({\rm mod}\ 3)\ \vee\ \beta\equiv 2\ ({\rm mod}\ 3)\}]$.

The latter notion of forcing is semiproper and $S$-closed,
being $\omega_1$-closed. Its introduction necessitates that we define
the $P_\alpha$'s with mixed support: Proceeding recursively, for
$\alpha\leq\kappa$ let $\overline{P}_\alpha$ consist of functions $p$
with domain $\alpha$ such that for each $\beta <\alpha,\ p(\beta )$ is
a $\overline{P}_\beta$-name such that $\vd_{\overline{P}_\beta}\
p(\beta )\in\dot R_\beta$, and 
$$
\eqalign{|{\rm supt}(p )&\cap \{
\beta <\alpha\ |\ \beta\equiv 0\ ({\rm mod}\ 3)
\ \vee\ \beta\equiv 1\ ({\rm mod}\ 3)\}|\leq\aleph_0\cr
|{\rm supt}(p )&\cap \{ \beta <\alpha\ |\ \beta\equiv 2\ ({\rm mod}\ 3)\}|\leq
\aleph_1.\cr}
$$
The following lemma will be a consequence of forthcoming iteration lemmas.
\bigskip
{\bf 3.1 Lemma.} {\it For each $\alpha\leq\kappa,\
\overline{P}_\alpha$ is $S$-closed and semiproper.
}
\bigskip
Assuming this lemma, the proof of Theorem C can be completed as
follows: For $\alpha\leq \kappa$, define
$$\eqalign{
\overline{P}_\alpha^\prime =\{ p\in\overline{P}_\kappa\
|\ &\forall\alpha\in {\rm supt}(p)(\alpha\equiv 1\ ({\rm mod}\ 3)\ \lor\ 
\alpha\equiv 2\ ({\rm mod}\ 3)) \cr
&\land\forall \beta\in {\rm
supt}(p)(p(\beta)\ \hbox{\rm is a}\ P^\prime_\beta{\rm -name})\}\;.
\cr}$$ 
($p(\beta )$ being a $P^\prime_\beta$-name, it only depends on $\dot G_\beta
\cap P^\prime_\beta$.) Then any forcing extension via $\overline{P}_\kappa
^\prime$ satisfies the conclusion of Theorem C:

Let $G^\prime$ by $\overline{P}_\kappa^\prime$-generic. It can be
checked that the analogue of 2.4 holds in the new situation. In
particular, there is a $G$ $\overline{P}_\kappa$-generic such that
$G\cap \overline{P}_\kappa^\prime = G^\prime$. With 3.1, the proof of
2.1 still works to show that in $V[G],\
\kappa =\omega_2$ and ${\rm NS}|\tilde S$ is $\kappa$-saturated. It then follows as
 before that in $V[G^\prime ]$, NS$|\tilde S$ is $\kappa$-saturated and 
${\cal P}(\omega_1)/{\rm NS}|\tilde S$ is a complete Boolean algebra.

By standard arguments $\overline{P}_\kappa^\prime$ has the
$\kappa$-c.c., and the introduction of the collapses $\dot R_\alpha$
for $\alpha\equiv 2$ (mod 3) implies that in $V[G^\prime ],\ \kappa
=2^{\omega_1}\leq\omega_3$. But by 3.7 below, $\omega_2$ is preserved so
that $\kappa =\omega_3$ in $V[G^\prime ]$ and so the proof is
complete.

\vskip 0.2in
The rest of this section is devoted to establishing 3.1 and the forthcoming 
3.7. We build on Shelah [Sh2][Sh4:~XIV] and refer to them for the more basic 
details about iterated forcing that are not provided in full.

For $\mu >\omega$, let ${\cal K}_\mu$ be the class of $\langle
Q,\leq_Q,\leq_Q^0\rangle$ such that:

(i) $\langle Q,\leq_Q\rangle$ is a semiproper, $S$-closed notion of
forcing, say with maximum element $\dot 1_Q$.

(ii) $\langle Q,\leq_Q^0\rangle$ is a poset so that: (a) $\leq_Q$
refines $\leq_Q^0$, i.e. if $p\leq_Q^0 q$, then $p\leq_Q q$; and (b)
$\leq_0$ is $\mu$-closed, i.e. if $\langle p_\alpha\ |\ \alpha
<\eta\rangle$ is $\leq_Q^0$-decreasing and $\eta <\mu$, then there is
a $p\in Q$ such that $p\leq_Q^0 p_\alpha$ for every $\alpha <\eta$.

We often suppress the subscripts $Q$ and furthermore identify $\langle
Q,\leq_Q,\leq_Q^0\rangle$ with its domain $Q$ when there is little
possibility of confusion. When we use forcing terminology for such a
member of ${\cal K}_\mu$, we are referring to the $\langle
Q,\leq_Q\rangle$ part.

Next, let ${\cal K}_\mu^\ast$ be the class of $\langle P_\alpha,\dot
Q_\alpha\ |\ \alpha <\gamma\rangle$ where for each $\alpha <\gamma$,
$P_\alpha$ is a notion of forcing, $\dot Q_\alpha$ is a
$P_\alpha$-name and $\vd_{P_\alpha}\
\dot Q_\alpha\in {\cal K}_\mu$ and recursively:

(i) $P_\alpha$ consists of functions $p$ with domain $\alpha$ so that
for each $\beta <\alpha,\ p(\beta )$ is a $P_\beta$-name such that
$\vd\ p(\beta )\in \dot Q_\beta$, and setting 
$${\rm supt}(p)=\{ \beta
<\alpha\ |\ \vd\ p(\beta )\not= \dot 1_{Q_\beta}\}$$ 
as before, 
$$|{\rm supt}(p)|\leq\aleph_1$$ 
and 
$$|\{ \beta <\alpha\ |\ \lnot (\vd\ p(\beta )\leq_{\dot Q_\beta}^0 \dot 
1_{Q_\beta})\}|\leq\aleph_0\;.
$$

(ii) The ordering on $P_\alpha$ is given by: 
$$\eqalign{
p\leq q\ \hbox{\it iff\/}\ &\forall \beta <\alpha (p|\beta\ \vd \ p(\beta ) 
\leq_{\dot Q_\beta} q(\beta ))\ \ \land\cr
&|\{ \beta\in {\rm supt}(p)\ |\ \lnot (p|\beta\ \vd\ p(\beta )\ 
\leq_{\dot Q_\beta}^0\ q(\beta ))\}|\leq\aleph_0\;.\cr}
$$
If $P_\gamma$ is defined by taking $\alpha =\gamma$ in the above, we
say that $P_\gamma$ is the {\it iteration} of $\langle P_\alpha,\dot
Q_\alpha\ |\ \alpha <\gamma\rangle$. The introduction of the second
partial order $\leq^0$ serves less to provide iteration lemmas of
potentially wide applicability than to provide a uniform approach to
3.1. For that result, $\leq_{\dot R_\alpha}^0$ will coincide with
$\leq_{\dot R_\alpha}$ when $\alpha\equiv 2$ (mod 3), i.e. when $\dot
R_\alpha$ is the Levy collapse col$(\omega_2,2^{\omega_2})$; and
$\leq_{\dot R_\alpha}^0$ will just be equality for $\alpha\equiv 0$
(mod 3) and $\alpha\equiv 1$ (mod 3). Note then that by how we defined
the $\overline{P}_\alpha$'s from the $\dot R_\alpha$'s,
$\langle\overline{P}_\alpha,\dot R_\alpha\ |\ \alpha <\kappa\rangle\in
{\cal K}_\mu^\ast$.

The usual iteration facts hold for members of ${\cal K}_\mu^\ast$. For
example, if $\langle P_\alpha,\dot Q_\alpha\ |\ \alpha <\gamma\rangle
\in {\cal K}_\mu^\ast$, for any $\beta <\alpha <\gamma,\ P_\beta \lessdot
P_\alpha$, i.~every maximal antichain of $P_\beta$ is a maximal
antichain of $P_\alpha$, and for the usual quotient poset $P_\alpha
/P_\beta$ such that $P_\alpha\cong P_\beta\ast P_\alpha /P_\beta$,
$P_\alpha /P_\beta$ is an iteration of a member of ${\cal
K}_\mu^\ast$. To establish 3.1, we must verify that iterations in
${\cal K}_\mu^\ast$ preserve $S$-closure and semiproperness.

The following lemma provides the main induction step for establishing
the preservation of $S$-closure:
\bigskip
{\bf 3.2 Lemma.} {\it Suppose that\/ $P_\gamma$ is the iteration of\/ 
${\cal P}= \langle P_\alpha,\dot Q_\alpha\ |\ \alpha <\gamma\rangle\in 
{\cal K}_\mu^\ast$ and $\delta <\eta\leq\gamma$. Then the following holds for
sufficiently large regular $\lambda$:

Assume that $N$ is countable with $\langle N,\in\rangle\prec
\langle H_\lambda,\in\rangle$, $N\cap\omega_1\in S,\ \{ {\cal
P},\delta,\eta\}\subseteq N$, $\langle p_n\ |\ n\in\omega\rangle$ is a
$P_\eta$-generic sequence for $N$, and $q\in P_\delta$ satisfies
$q\leq p_n|\delta$ for every $n\in\omega$.  Then there is a $q^{+}\in
P_\eta$ such that $q^{+}|\delta =q$ and $q\leq p_n$ for every
$n\in\omega$.}
\medskip
{\bf Proof}. By adjusting names, we can assume for convenience in what
follows that for each $n\in\omega$ and $\alpha <\eta$:

(i) $p_{n+1}|\alpha\ \vd\ p_{n+1}(\alpha )\leq^0p_n(\alpha )\ \ {\it iff}\ \ 
\vd\ p_{n+1}(\alpha )\leq^0 p_n(\alpha )$, and

(ii) $\vd\ p_{n+1}(\alpha )\leq p_n(\alpha )$.

We now define a function $q^{+}$ with domain $\eta$ as follows: Fix a
well-ordering $W$ of a sufficiently large $V_\rho$. Set $q^{+}|\delta
=q$. For $\delta\leq\alpha <\eta$, having defined $q^{+}|\alpha$ for
$\alpha <\eta$ so that recursively $q^{+}|\alpha\in P_\alpha$, define
$q^{+}(\alpha )$ as follows:

(a) If for some $k\in\omega$, $q^{+}|\alpha\ \vd\ ``\langle
p_n(\alpha )\ | \ k<n<\omega\rangle$ is $\leq^0$-descending in $\dot
Q_\alpha - \{\dot 1_{Q_\alpha}\}$'', then since by definitions of
${\cal K}_\mu$ and ${\cal K}_\mu^\ast$, $\vd_{P_\alpha}\ ``\leq^0_{Q_\alpha}$ 
is $\mu$-closed'' and
$\mu >\omega$, there is a $P_\alpha$-name $\tau$ so that
$q^{+}|\alpha\ \vd\ ``\tau$ is a $\leq^0$-lower bound for $\langle
p_n(\alpha )\ |\ k<n<\omega\rangle$''.  Let $q^{+}(\alpha )$ be the
$W$-least such $\tau$. Else:

(b) If for some $k\in\omega$ and $P_\alpha$-name $\tau$,
$q^{+}|\alpha\ \vd\ ``\langle p_n(\alpha )\ |\ k<n<\omega\rangle$ is
$\leq_{Q_\alpha}$-descending in $\dot Q_\alpha - \{ 1_{Q_\alpha}\}$
with $\tau$ a $\leq_{Q_\alpha}$-lower bound'', then let $q^{+}(\alpha )$ 
be the $W$-least such $\tau$. Otherwise:

(c) Set $q^{+}(\alpha )=\dot 1_{Q_\alpha}$.

This definition perpetuates $q^{+}|\alpha\in P_\alpha$ for every
$\alpha\leq\eta$: Clearly, 
$$
|{\rm supt}(q^{+}|\alpha)|\leq |{\textstyle \bigcup}_n{\rm supt}
(p_n|\alpha )|\leq \aleph_1\;.
$$ 
Also,
conditions on $\leq^0$ and (i) above imply that there is a countable
set $E$ such that for $\beta\in\alpha - E$, $\vd\ ``p_{n+1}(\beta)
\leq^0 p_n(\beta )\leq^0 \dot 1_{Q_\beta}$'' for every $n\in\omega$. For such
$\beta,\ q^{+}(\beta )$ was defined either through clause (a), or if
not, neither through clause (b) as $\leq$ refines $\leq^0$, but
through clause (c).  But both (a) and (c) lead to $\vd\ q^{+}(\beta
)\leq^0\dot 1_{Q_\zeta}$, and so 
$$\{ \beta <\alpha\ |\ \lnot (\vd\
q^{+}(\beta )\leq^0\ \dot 1_{Q_\beta})\}\subseteq E,$$
confirming that $q^{+}|\alpha\in P_\alpha$.

We next establish by induction on $\alpha\leq\eta$ that

(1) for every $n\in\omega,\ q^{+}|\alpha\ \vd\ q^{+}(\alpha )\leq
p_n(\alpha )$. 

(2) If $\alpha\not\in N - \delta$, then for every $n\in\omega,\
q^{+}|\alpha\ \vd\ q^{+}(\alpha )\leq^0p_n(\alpha )$.

As $N$ is countable, this suffices to verify that $q^{+}\leq p_n$ for
every $n\in\omega$ as desired. For $\alpha <\delta$, the results are
immediate; what remains splits into two cases:

{\bf Case 1.} $\delta\leq\alpha <\eta$ and $\alpha\in N$. Then
$P_\alpha,\dot Q_\alpha\in N$, and we have $q|\alpha\leq p_n|\alpha$
for every $n\in \omega$ by induction. Since $\langle p_n\ |\
n\in\omega\rangle$ is a $P_\eta$-generic sequence for $N$, it follows
first that $\langle p_n|\alpha\ | \ n\in\omega\rangle$ is a
$P_\alpha$-generic sequence for $N$, and second that $q^{+}|\alpha\
\vd\ ``\langle p_n(\alpha )\ |\ n\in\omega\rangle$ is a $\dot
Q_\alpha$-generic sequence for $\check N[\dot G_\alpha ]$''. But by
definitions of ${\cal K}_\mu$ and ${\cal K}_\mu^\ast$,
$\vd_{P_\alpha}\ ``\dot Q_\alpha$ is $S$-closed'', and since $N\cap
\omega_1\in S$, (a) or (b) of the definition of $q^{+}$ applies at
$\alpha$ and (1) follows. (2) holds vacuously.

{\bf Case 2.} $\delta\leq\alpha <\eta$ and $\alpha\not\in N$. The
$\leq^0$ conditions on the $p_n$'s imply that for each $n\in\omega$,
$$E_n=\{ \beta <\eta\ |\ \lnot (\vd\ p_{n+1}(\beta )\leq^0 p_n(\beta
)\leq^0\dot 1_{Q_\beta})\}$$ 
is countable, and clearly $E_n\in N$, so that $E_n$ is countable in
$N$. Hence, $\bigcup_nE_n\subseteq N$, so that $\alpha\not\in
\bigcup_nE_n$. Hence, $q^{+}(\alpha )$ was defined either through
clause (a), or if not, neither through clause (b) as $\leq$ refines
$\leq^0$, but through clause (c). Both (1) and (2) now follow in this
case also. \hfill $\dashv$

The following preservation result now follows in straightforward fashion:
\bigskip
{\bf 3.3 Proposition.} {\it Suppose that $P_\gamma$ is the iteration of 
$\langle P_\alpha,\dot Q_\alpha\ |\ \alpha <\gamma\rangle\in {\cal
K}_\mu^\ast$ and $\delta <\eta\leq\gamma$. Then $P_\eta /P_\delta$ is
$S$-closed.} 
\bigskip
The following lemma provides the main induction step for establishing
the preservation of semiproperness. Again, for a notion of
forcing $P$, $\dot G_P$ denotes the canonical $P$-name for its
generic object. For a set $M$ and $q\in P$, to say that $q$ {\it
decides} $\dot G_P\cap \check M$ means of course that there is a set
$A\in V$ such that $q\ \vd\ \check A=\dot G_P\cap \check M$. This
happens for example if $M$ is countable and $q$ is a lower bound to a
$P$-generic sequence for $M$, with $A=\{ r\in P\cap M\ |\ q\leq r\}$.
Finally, for $\lambda >\omega$ regular and $\langle N,\in\rangle
\prec \langle H_\lambda,\in\rangle$ with $P\in N$, if $G$ is
$P$-generic over $V$, then let $N[G]$ be the set of interpretations
$\{ \tau [G]\ |\ \tau\in N\ \hbox{\rm is a $P$-name}\}$. By Shelah
[Sh2: p.88], 
$$\langle N[G],\in\rangle\ \prec\ \langle
H_\lambda^{V[G]},\in\rangle.\leqno (\ast )$$
We let $N\dot [G_P]$ be a canonical $P$-name for $N[G]$.
\bigskip
{\bf 3.4 Lemma.} {\it Suppose that $P_\gamma$ is the iteration of
${\cal P}=\langle P_\alpha,\dot Q_\alpha\ |\ \alpha <\gamma\rangle\in
{\cal K}_\mu^\ast$ and $\delta <\eta\leq\gamma$. Then the following
holds for sufficiently large regular $\lambda$:

Assume that $N$ and $M$ are countable with $\langle N,\in\rangle
\prec \langle M,\in\rangle \prec \langle H_\lambda,\in\rangle$, $\{
{\cal P}, \delta,\eta\}\subseteq N$, $N\in M$, and $M\cap\omega_1\in
S$. Assume also that $q\in P_\delta$ is $\langle
N,P_\delta\rangle$-semigeneric, $q$ is $\langle M,
P_\delta\rangle$-generic and decides $\dot G_{P_\delta} \cap\check
M$, and $p\in P_\eta\cap N$ with $q\leq p|\delta$. 

Then there is a $q^{+}\in P_\eta$ with $q^{+}\leq p$ and $q^{+}|\delta =q$ 
such that $q^{+}$ is $\langle N,P_\eta\rangle$-semigeneric, and $q^{+}$ is 
$\langle M,P_\eta\rangle$-generic and decides $\dot G_{P_\eta}\cap \check M$.}
\medskip
{\bf Proof}. We establish this by induction on $\eta$, for all $\delta,N,M,p$,
and $q$.

For $\eta$ a successor, we can clearly assume that $\eta =\delta +1$. By 
definitions of $K_\mu$ and ${\cal K}_\mu^\ast$, $\vd_{P_\delta}\
``\dot Q_\delta$ is semiproper'', and so with $(\ast )$ just before
3.1 in mind, we have 
$$\vd_{P_\delta}\ \exists r\in\dot Q_\delta (r\ \hbox{\rm is}\ \langle N\dot [
G_{P_\delta}],P_\delta\rangle\hbox{\rm -semigeneric}\ \land\ r(\delta )\leq
 p(\delta )).\leqno (\ast\ast )$$
Since $N\in M$ and $q$ is $\langle M,P_\delta\rangle$-generic and
decides $\dot G_{P_\delta}\cap\check M$, by $(\ast )$ applied
syntactically there is a $q^\prime\in P_\delta\cap M$ with $q\leq q^\prime$
and a $P_\delta$-name $r^\prime\in M$ such that $q^\prime$ forces the
assertion of $(\ast\ast )$ with $r^\prime$.  Since $M\cap\omega_1\in
S$, $\vd_{P_\delta}\ ``\dot Q_\delta$ is $S$-closed'' by the definitions,
and with $(\ast )$ in mind for $M$, it is now straightforward to find
a $q^{+}\in P_\eta$ as desired, noting that in the sense of
$\vd_{P_\delta}$, any lower bound to a $\dot Q_\delta$-generic
sequence for $\check M$ decides $\dot G_{Q_\delta}\cap\check M$.

For a $\eta$ a limit, let $\langle\beta_n\ |\ n\in\omega\rangle$ enumerate $M
\cap \{ \beta\ |\ \delta\leq\beta <\eta\}$ with $\beta_0=\delta$. Let
$\{\dot{\tau}^l_n: n<\omega\}$ enumerate the $P_\eta$-names of
ordinals belonging to $N$ if $l=1$ and $M$ if $l=2$. Define 
$\alpha_n,q_n,p_n$, and $N_n$ by induction on $n\in\omega$ so that henceforth 
writing $\dot G_n$ for $\dot G_{P_{\alpha_n}}$, we have:

(a) $\alpha_0=\delta,\ q_0=q$, $p_0=p$ and generally $q_n\in
P_{\alpha_n}$ with $q_n \leq p|\alpha_n$ and $q_{n+1}|\alpha_n=q_n$,
$p_n\in N[\dot{G}_n]\cap P_{\alpha_n}$ ($\subseteq M\cap P_\alpha$),
$p_{n+1}\leq p_n$, and $p_{n+1}$ decides a value for $\dot{\tau}_n$.

(b) $q_n$ is $\langle N,P_{\alpha_n}\rangle$-semigeneric, and $q_n$ is 
$\langle M,P_{\alpha_n}\rangle$-generic and decides $\dot G_n\cap\check M$.

(c) $q_n\ \vd\ \alpha_{n+1}= \max\{ \{ \beta_k\ |\ k\leq n\}\cap N\dot [ G_n]
\}$.

\noindent (Note for (c) that if $q_n$ decides $\dot G_n \cap\check M$ by (b),
then $N\in M$ implies that $q_n$ decides $\dot G_n \cap\check N$ and
hence $N\dot [ G_n]$.)

The case $n=0$ follows from our initial assumption. Suppose now that
$q_n$ has already been defined. By $(\ast )$ just before 3.4 applied
syntactically, $\vd\ ``N\dot [ G_n]\prec M\dot [ G_n]\prec \langle 
H_\lambda^{V[\dot G_n]},\in\rangle''$. Moreover, since
$q_n$ is $\langle M, P_{\alpha_n}\rangle$-generic, we have $q_n\ \vd\
M\dot [ G_n]\cap\omega_1= \check M\cap\omega_1\in \check S$. Hence,
with $\alpha_{n+1}$ as stipulated by (c) it is straightforward to
apply the induction hypothesis in the sense of $q_n\ \vd$ and then to
find an appropriate $q_{n+1}\in P_{\alpha_{n+1}}$ as desired. There is no
problem in defining $p_{n+1}$. 

We can now define a $q^{+}\in P_\eta$ so that ${\rm
supt}(q^{+})=\bigcup_n {\rm supt}(q_n)$, and for any $\beta$ in this
set, $q^{+}(\beta )=q_n(\beta )$ for some (any) $n$ such that
$\beta\in {\rm supt}(q_n)$. As in [Sh2][Sh4], $q \leq p_n$ so $q^+$ is 
$\langle N,P_\alpha\rangle$-semigeneric.
$q^+$ is not necessarily $(M,P_\alpha)$-generic, but its existence
shows that there is $q^\prime\in P_\eta\cap M$, $q^\prime |
\alpha\leq q^+$, $q^\prime$ is $\langle N,P_\alpha\rangle$-semigeneric. Now 
we can find a $q^+$ really as required.\hfill $\dashv$
\bigskip
The following preservation result is now clear, since for $\langle N,\in\rangle
\prec \langle H_\lambda,\in\rangle$ as in 3.4, we can always find 
a countable $M$ such that $N\in M,\ \langle M,\in\rangle \prec \langle 
H_\lambda,\in\rangle$, and $M\cap\omega_1\in S$ by the stationariness of $S$.
\bigskip
{\bf 3.5 Proposition.} {\it Suppose that $P_\gamma$ is the iteration of 
$\langle P_\alpha,\dot Q_\alpha\ |\ \alpha <\gamma\rangle\in {\cal K}_\mu^\ast$
 and $\delta <\eta\leq\gamma$. Then $P_\eta /P_\delta$ is semiproper.}
\bigskip
{\bf Proof of 3.1.} By 3.3 and 3.5 we get $S$-completeness and 
semiproperness.\hfill $\dashv$ 
\bigskip
Finally, we establish the preservation of $\omega_2$ in a special
case; a similar result appears in Shelah \hfil\break
[Sh2:~VIII,\S1].
\bigskip

\bigskip
{\bf 3.6 Proposition.} {\it Suppose that {\rm CH} holds, and $P$ is the
 iteration of\/ ${\cal P}=\langle P_\alpha,\dot Q_\alpha\ |\ \alpha <\gamma\rangle\in 
{\cal K}_{\omega_2}^\ast$ where for each $\alpha <\gamma$, $\vd_{P_\alpha}\ ``
{\leq}_{\dot Q_\alpha}^0= {\leq}_{\dot Q_\alpha}$ and $Q_\alpha$ satisfies
$\aleph_2$-pic'' (see [Sh2:~VII, Def.2.]) or\/
$\vd_{P_\alpha}\ ``\leq_{\dot Q_\alpha}^0$ is the equality relation on
$\dot Q_\alpha$''. Then forcing with $P$ preserves $\omega_2$.}
\medskip
{\bf Proof.} Assume to the contrary that for some $p\in P$ and $P$-name $\tau$,
$p\ \vd\ ``\tau\colon\check\omega_1\rightarrow\check\omega_2$ is a 
bijection''. Taking a regular $\lambda$ sufficiently large, we proceed by 
induction on $\xi <\omega_2$ to define $p_\xi\in P$ so that $\xi <\zeta <
\omega_2$ implies that for every $\alpha <\gamma,\ \vd\ p_\zeta (\alpha )\leq^0
p_\xi (\alpha )$, and accompanying $q_\xi\leq p_\xi,\ N_\xi$, and $A_\xi$ as
follows:

Set $p_0=p$. At limits $\zeta <\omega_2$ with $p_\xi$ having been defined for
$\xi <\zeta$, by the $\omega_2$-closure of the $\leq_{\dot Q_\alpha}^0$'s,
let $p_\zeta\in P$ be such that for every $\xi <\zeta$ and $\alpha <\gamma$, 
$\vd\ p_\zeta (\alpha )\leq^0 p_\xi (\alpha )$.

To handle the successor stage, suppose that $p_\xi$ is given. First let $N_\xi$
be countable with $\langle N_\xi,\in\rangle\prec \langle H_\lambda,\in
\rangle$, $\{ {\cal P},\xi, p_\xi\}\subseteq N_\xi$, and $N_\xi\cap\omega_1\in
S$. Then define $q_\xi\leq p_\xi$ by first choosing a $P$-generic sequence for
$N_\xi$ starting with $p_\xi$ and then using the clauses (a),(b), and (c) as
in the proof of 3.2. Hence $q_\xi$ is $\langle N_\xi, P\rangle$-generic and
setting $A_\xi =\{ r\in P\cap N_\xi\ |\ q_\xi\leq r\}$, we have $q_\xi\ \vd\ 
\check A_\xi =\dot G_P\cap\check N_\xi$. We can assume that $\vd\ q_\xi(\alpha 
)\ \leq^0\ p_\xi (\alpha )$\ {\it iff\/}\ $q_\xi |\alpha\ \vd\ q_\xi 
(\alpha ) \leq^0 p_\xi (\alpha )$; similarly for $\leq$; and also that each 
$q_\xi (\alpha )$ depends only on $p_\xi (\alpha )$ and not on $\xi$ 
(by carrying out
the entire procedure canonically according to some well-ordering of a 
sufficiently large $V_\rho$). Finally, define $p_{\xi +1}$ as follows:
$$p_{\xi +1}(\alpha )=\cases{q_\xi (\alpha )&if $\vd\ q_\xi (\alpha )\leq^0
p_\xi (\alpha )$,\cr
\hfil\cr
p_\xi (\alpha )&otherwise.\cr}$$
\indent Proceeding with the proof, a straightforward $\Delta$-system argument
with CH shows that there are $\xi_0 <\xi_1 <\omega_2$ such that: $N_{\xi_0}
\cap\omega_1=N_{\xi_1}\cap\omega_1$, and there is an isomorphism $h\colon 
\langle N_{\xi_0},\in,A_{\xi_0}\rangle\rightarrow\langle N_{\xi_1},\in,
A_{\xi_1}\rangle$ with $h|(N_{\xi_0}\cap N_{\xi_1})$ the identity, $h({\cal P})
={\cal P},\ h(\xi_0)=\xi_1$, and $h(p_{\xi_0})=p_{\xi_1}$. By assumption, there
is an $\eta\in N_{\xi_0}\cap\omega_1$ such that $q_{\xi_0}\ \vd\ \dot\tau 
(\eta )=\xi_0$. Hence, $\exists r\in A_{\xi_0}(r\ \vd\ \dot\tau (\eta )=
\xi_0)$. Applying $h$, we have $\exists r\in A_{\xi_1}(r\ \vd\ \dot\tau (\eta )
=\xi_1)$ since $h$ is the identity on $N_{\xi_0}\cap\omega_1$. Consequently,
$q_{\xi_1}\ \vd\ \dot\tau (\eta )=\xi_1$. However, it is not difficult to 
check that $q_{\xi_0}$ and $q_{\xi_1}$ are compatible, reaching
a contradiction:

The definitions of $q_{\xi_0}$ and $q_{\xi_1}$ as in the proof of 3.2 show that
the countable sets
$$
E_i=\{ \alpha\in {\rm supt}(q_{\xi_i})\ |\ \lnot (\vd\ q_{\xi_i}(\alpha ) 
\leq^0 p_{\xi_i}(\alpha ))\}
$$
for $i<2$ are such that $E_i\in N_{\xi_i}$ and $E_i$ is countable in 
$N_{\xi_i}$. We now argue that $\vd\ ``q_{\xi_0}(\alpha )$ and $q_{\xi_1}
(\alpha )$ are compatible'' by cases, depending on whether 
$\alpha\in\gamma -
E_0,\ \alpha\in E_0 - E_1$, or $\alpha\in E_0\cap E_1$: If $\alpha\in\gamma
- E_0$, then $\vd\ q_{\xi_1}(\alpha )\leq p_{\xi_1}(\alpha )\leq^0 
p_{\xi_0+1}(\alpha )\leq^0 q_{\xi_0}(\alpha )$. If $\alpha\in E_0 - E_1$, 
then $\vd\ q_{\xi_0}(\alpha )\leq p_{\xi_0}(\alpha )$ and $\alpha\in E_0$ 
implies that we are in the case $\vd_{P_\alpha} ``\leq_{\dot Q_\alpha}^0$ is
equality'' of our assumption about $P$. Consequently, $\alpha\not\in E_1$ 
implies that $\vd\ q_{\xi_1}(\alpha )=p_{\xi_1}(\alpha )=p_{\xi_0}(\alpha )$,
and so we have $\vd\ q_{\xi_0}(\alpha )\leq q_{\xi_1}(\alpha )$. Finally, if
$\alpha\in E_0\cap E_1$, then $\alpha\in N_{\xi_0}\cap N_{\xi_1}$ 
so the first case in the proposition occurs and we apply ``$\aleph_2$-pic''.
This completes the proof.\hfill $\dashv$

\bigskip
In 3.6 we can combine the two possibilities to one as
implicit in the proof.
\bigskip

{\bf 3.7 Proposition.} {\it $P_\kappa^\prime$ (from the proof of Theorem C) 
preserves $\omega_2$.}
\medskip
{\bf Proof.} Let $\langle P^1_\alpha,\dot{Q}^1_\alpha:
\alpha<\kappa\rangle\in {\cal K}_\mu$ be as above except that for 
$\alpha\equiv 0$ (mod 3), $\dot Q^1_\alpha$ is the trivial forcing (and if 
$\alpha\not\equiv 0$ (mod 3), $\dot{Q}^1_\alpha=\dot{R}_\alpha$). Essentially,
$P^1_\alpha=P_\alpha^\prime$, hence it suffices to
prove that $P^1_\alpha$ preserve $\aleph_2$. Now the
assumption of 3.6 clearly holds for
$\langle P^1_\alpha,\dot{Q}^1_\alpha: \alpha<\kappa\rangle$ (for
$\aleph_2$-pic --- see the discussion of $Q_S$ for more). \hfill $\dashv$
\bigskip
This finally completes the proof of Theorem C. 

Instead of $2^{\omega_1}=
\omega_3$, for any regular $\nu$ such that $\omega_2\leq\nu <\kappa$, it is
possible to arrange $2^{\omega_1}=\nu^{+}=\kappa$ using ${\cal K}_\nu^\ast$
with $|{\rm supt}(p)|\leq\nu$ in place of $|{\rm supt}(p)|\leq\omega_1$ in
its definition.

\vfill\eject
\centerline{\bf References}
\medskip
\item{[B]} Baumgartner, James, Iterated Forcing, in: Mathias (ed.) {\it Surveys
 in Set Theory}, London Mathematical Society Lecture Note Series \#87 
(Cambridge University Press, 1983), 1-59.
\medskip
\item{[BT]} Baumgartner, James and Alan Taylor, Saturation properties of ideals
 in generic extensions, II, {\it Trans. Amer. Math. Soc.} {\bf 27} (1982),
 587-609.
\medskip
\item{[FMS]} Foreman, Matthew, Menachem Magidor, and Saharon Shelah, Martin's
 Maximum, saturated ideals and non-regular ultrafilters. Part I, {\it Annals
 Math.} {\bf 127} (1988), 1-47.
\medskip
\item{[GS]} Gitik, Moti and Saharon Shelah, Cardinal preserving ideals, to 
appear in the {\it Israel Jour. Math.}
\medskip
\item{[J]} Jech, Thomas, {\it Set Theory} (New York: Academic Press, 1978).
\medskip
\item{[Ke]} Ketonen, Jussi, Some combinatorial principles, {\it Trans. Amer.
 Math. Soc.} {\bf 188} (1974), 387-394.
\medskip
\item{[Ku]} Kunen, Kenneth, Some applications of iterated ultrapowers in set
 theory, {\it Ann. Math. Logic} {\bf 1} (1970), 179-227.
\medskip
\item{[P]} Pierce, R.S., Distributivity in Boolean algebras, {\it Proc. Amer. 
Math. Soc.} {\bf 7} (1957), 983-992.
\medskip
\item{[Sh1]} Shelah, Saharon, A weak generalization of MA to higher cardinals,
 {\it Israel Jour. Math.} {\bf 30} (1978), 297-306.
\medskip
\item{[Sh2]} Shelah, Saharon, {\it Proper Forcing}, Lecture Notes in 
Mathematics \#940 (Berlin: Springer-Verlag, 1986).
\medskip
\item{[Sh3]} Shelah, Saharon, Iterated forcing and normal ideals on $\omega_1$,
 {\it Israel Jour. Math.} {\bf 60} (1987), 345-380.
\medskip
\item{[Sh4]} Shelah, Saharon, {\it Proper and Improper Forcing} (Berlin: 
Springer-Verlag), to appear. Revised and expanded version of [Sh2].
\medskip
\item{[Si]} Sikorski, Roman, On an unsolved problem from the theory of Boolean
 algebras, {\it Coll. Math.} {\bf 2} (1949), 27-29.
\medskip
\item{[ST]} Smith, Edgar and Alfred Tarski, Higher degrees of distributivity
 and completeness in Boolean algebras, {\it Trans. Amer. Math. Soc.} {\bf 84}
 (1957), 230-257.
\medskip
\item{[So]} Solovay, Robert, Real-valued measurable cardinals, in: Scott, Dana
 (ed.) {\it Axiomatic Set Theory}, Proceedings of Symposia in Pure Mathematics
 vol. 13, part 1 (Providence; A.M.S., 1971), 397-428.
\bye